\documentclass[11pt,a4paper,fullpage]{article}
\usepackage{amsmath,amsthm}
\usepackage{amssymb}

\newtheorem{theorem}{Theorem}

\title{{\bf Dedekind sums and class numbers of imaginary abelian number fields}}
\author{
St\'ephane R. LOUBOUTIN\\
Aix Marseille Universit\'e, CNRS, Centrale Marseille, I2M,\\ 
Marseille, France\\
stephane.louboutin@univ-amu.fr}
\date{\today}

\begin{document}
\bibliographystyle{alpha}
\maketitle


\begin{abstract}
As a consequence of their work, 
Bruce C. Berndt and Ronald J. Evans in 1977 and Larry Joel Goldstein and Michael Razar in 1976 
obtained a formula for the square of the class number of an imaginary quadratic number field 
in terms of Dedekind sums. 
Here we focus only on this formula and give a new short and simple proof of it. 
We finally express the relative class numbers of imaginary abelian number fields in terms of Dedekind sums.
 \end{abstract}

For $c,d\in {\mathbb Z}$ and $d>1$ such that $\gcd (c,d)=1$,
the {\it Dedekind sum} is defined by
\begin{equation}\label{defs}
s(c,d)
:={1\over 4d}\sum_{n=1}^{\vert d\vert -1}\cot\left ({\pi n\over d}\right )\cot\left ({\pi nc\over d}\right )
\end{equation}
(see \cite[Chapter 3, Exercise 11]{Apo} or \cite[(26)]{RG}).
It depends only on $c$ modulo $d$.

Let $3<p\equiv 3\pmod 4$ be a prime rational integer. 
The class number formula applied to the imaginary quadratic field 
$K ={\mathbb Q}(\sqrt{-p})$ of discriminant $-p$ yields 
$$h_K 
={\sqrt p\over\pi}L(1,\chi_K),$$
where $n\mapsto\chi_K(n) =\left (\frac{n}{p}\right )\in\{-1,0,+1\}$ is the odd an primitive Dirichlet character modulo $p$ 
associated with $K$ given by the Legendre symbol
and $L(s,\chi_K)$ is the Dirichlet $L$-function associated with the Dirichlat character $\chi_K$.
By \cite[Proposition 1]{LouCMB36/37} we have 
$$L(1,\chi_K)
=\frac{\pi}{2p}\sum_{n=1}^{p-1}\chi_K(n)\cot\left (\frac{\pi n}{p}\right ).$$
Consequently, we have 
\begin{eqnarray*}
h_K^2
=\frac{p}{\pi^2}L(1,\chi_K)^2
&=&\frac{\pi}{4p}\sum_{a=1}^{p-1}\sum_{b=1}^{p-1}\chi_K(ab)\cot\left (\frac{\pi a}{p}\right )\cot\left (\frac{\pi b}{p}\right )\\
&=&\frac{\pi}{4p}\sum_{a=1}^{p-1}\sum_{b=1}^{p-1}\chi_K(b)\cot\left (\frac{\pi a}{p}\right )\cot\left (\frac{\pi ab}{p}\right )
=\sum_{b=1}^{p-1}\chi_K(b)s(c,p),
\end{eqnarray*}
by the bijective change of variables $b\mapsto ab$ modulo $p$ 
and by noticing that $\chi_K(a)\chi_K(ab) =\chi_K(a)^2\chi_K(b)=\chi_K(b)$ for $\gcd (b,p)=1$.
We have finally proved a special case of \cite[Corollary 9]{BE} for prime discriminants. 
Since $S:=\sum_{c=1}^{p-1}s(c,p) =\sum_{c=1}^{p-1}s(p-c,p) =-S$, 
we have $S=0$ and 
$$h_K^2
=\sum_{b=1}^{p-1}(1+\chi_K(b))s(c,p)
=2\sum_{b=1\atop \chi_K(b) =+1}^{p-1}s(c,p)
=2S(H_2,p),$$ 
with the notation of \cite[Proposition 2]{LouBKMS56}.
By \cite[Theorem 6]{LouBKMS56} 
we recover the well known fact that the class number 
of the imaginary quadratic number field ${\mathbb Q}(\sqrt{-p})$ is odd for $3<p\equiv 3\pmod 4$ a prime integer.
We now give a short proof of \cite[Corollary 9]{BE} and \cite[p. 358]{GR}: 

\begin{theorem}
Let $-D<0$ be the discriminant of an imaginary quadratic field $K$. 
Let $w_K\in\{2,4,6\}$ and $h_K$ be the number of complex roots of unity contained in $K$ 
and the class number of $K$.
Let $\chi_K$ be the real, odd and primitive character modulo $D$ associated with $K$. 
Then 
$$h_K
=\frac{w_K^2}{4}
\sum_{n=1\atop\gcd (n,D)=1}^{D-1}\chi_K(n)s(n,D).$$
\end{theorem}

\begin{proof}
Recalling the class number formula 
$$h_K
=\frac{w_K\sqrt D}{2\pi}L(1,\chi_K),$$ e.g. see \cite[p. 38]{Was}, 
the result follows from the following Theorem \ref{L1X}.
\end{proof}

\noindent\frame{\vbox{
\begin{theorem}\label{L1X}
Let $\chi$ be an odd and primitive character modulo $f>2$. 
Then 
\begin{equation}\label{L1XDedekind}
\vert L(1,\chi)\vert^2
=\frac{\pi^2}{f}
\sum_{n=1\atop\gcd (n,f)=1}^{f-1}\chi (n)s(n,f).
\end{equation}
\end{theorem}
}}

\begin{proof}
Using \cite[Proposition 1]{LouCMB36/37} we get 
$$\vert L(1,\chi)\vert^2
=\frac{\pi^2}{4f}\sum_{a=1\atop\gcd (a,f)=1}^{f-1}\sum_{b=1\atop\gcd (b,f)=1}^{f-1}
\overline{\chi (a)}\chi(b)\cot\left (\frac{\pi a}{f}\right )\cot\left (\frac{\pi b}{f}\right ).$$
Making the bijective change of variables $b\mapsto ab\mod f$ we get 
$$\vert L(1,\chi)\vert^2
=\frac{\pi^2}{4f}\sum_{b=1\atop\gcd (b,f)=1}^{f-1}\chi(b)
\sum_{a=1\atop\gcd (a,f)=1}^{f-1}\cot\left (\frac{\pi a}{f}\right )\cot\left (\frac{\pi ab}{f}\right ).$$
Noticing that 
$$\sum_{\delta\mid a\atop \delta\mid f}\mu(\delta)
=\begin{cases}
1&\hbox{if $\gcd (a,f)=1$}\\
0&\hbox{if $\gcd (a,f)>1$}
\end{cases}$$
we obtain
$$\vert L(1,\chi)\vert^2
=\frac{\pi^2}{4f}\sum_{\delta\mid f}\frac{\mu(\delta)}{\delta}\sum_{b=1\atop\gcd (b,f)=1}^{f-1}\chi_K(b)s(b,f/\delta).$$
To arrive at the desired result it suffices to prove that for any $d<f$ dividing $f$ 
and for any primitive modulo $f$ Dirichlet character $\chi$ we have $S(\chi,f,d)=0$, where 
$$S(\chi,f,d)
:=\sum_{b\in ({\mathbb Z}/f{\mathbb Z})^*}\chi(b)s(b,d).$$
To this end we consider the canonical surjective morphism 
$\phi:({\mathbb Z}/f{\mathbb Z})^*\rightarrow ({\mathbb Z}/d{\mathbb Z})^*$ 
and notice that $s(b,d)$ depends only on $b$ modulo $d$. 
We obtain 
$$S(\chi,f,d)
=\sum_{B\in ({\mathbb Z}/d{\mathbb Z})^*}\left (\sum_{\phi(b)=B}\chi (b)\right )s(B,d)$$ 
and we explain why each sum $\sum_{\phi(b)=B}\chi (b)$ is equal to $0$.
Indeed, for a given $B\in ({\mathbb Z}/d{\mathbb Z})^*$, 
we pick up any $b_0\in ({\mathbb Z}/f{\mathbb Z})^*$ 
for which $\phi (b_0)=B$. 
Then $\phi(b)=B$ if and only if $b =b_0x$ with $x\in\ker\phi$. 
Hence, 
$\sum_{\phi(b)=B}\chi (b)
=\chi(b_0)\sum_{x\in\ker\phi}\chi (x)$. 
Finally, since $\chi$ is primitive modulo $f$, there exists $x_0\in\ker\phi$ with $\chi (x_0)\neq 1$ 
(otherwise $\chi$ would be induced by a character modulo $d$). 
Therefore the restriction of $\chi$ to the multiplicative group $\ker\phi$ is a non trivial character 
and it follows that $\sum_{x\in\ker\phi}\chi (x)=0$ 
(make the change of variables $x\mapsto x_0x$ in this sum).
\end{proof}

Since $\chi (n^*)=\overline{\chi(n)}$ and $s(n^*,f) =s(n,f)$, 
where $n^*$ is the inverse in the multiplicative group $({\mathbb Z}/f{\mathbb Z})^*$ of $n$ with $\gcd(n,f)=1$, 
the bijective change of variables $n\mapsto n^*$ shows that the right hand side of \eqref{L1XDedekind} 
being equal to its complex conjugate is indeed a real number. 
Notice also that the oddness of $\chi$ gives $\chi(f-n)s(f-n,f) =\chi(n)s(n,f)$ 
and
\begin{equation}\label{L1XDedekindbis}
\vert L(1,\chi)\vert^2
=\frac{2\pi^2}{f}
\sum_{n=1\atop\gcd (n,f)=1}^{f/2-1}\Re(\chi (n))s(n,f).
\end{equation}
Let $Q_K\in\{1,2\}$, $w_K$ and $h_K^-$ be the Hasse unit index, the number of complex roots of unity in $K$ 
and class number of an imaginary abelian number field $K$. 
Let $X_K^-$ be the set of primitive and odd Dirichlet characters associated with $K$.
The class number formula and the conductor-discriminant formula (see \cite[Theorem 3.11]{Was}) give 
$$h_K^-
=Q_Kw_K\prod_{\chi\in X_K^-}\frac{\sqrt{f_\chi}}{2\pi}L(1,\chi),$$
where $f_\chi$ denote the conductor of the primitive character $\chi$. 
Now, instead of using \cite[Theorem 4.17]{Was} to compute $h_K^-$ 
we can use this relative class number formula and \eqref{L1XDedekindbis}. 
For example, if $K$ is an imaginary cyclic quartic field of conductor $f$ 
and $\chi$ is any one of the two conjugate odd and primitive modulo $f$ quartic characters associated with $K$ 
we have $Q_K=1$ and 
$$h_K^-
=w_K\frac{f^2}{4\pi^2}\vert L(1,\chi)\vert^2
=\frac{w_K}{2}\sum_{n=1\atop\gcd (n,f)=1}^{f/2-1}\Re(\chi (n))s(n,f).$$
For example, for $K={\mathbb Q}(\zeta_5)$ we have $Q_K=1$ by \cite[Corollary 4.13]{Was}, $w_K=10$ and 
since $s(1,5)=1/5$, $s(2,5)=s(3,5)=0$ and $s(4,5)=-1/5$ we get
$$h_K^-
=5\sum_{n=1}^2\Re(\chi (n))s(n,5)
=5\Bigl (1\cdot s(1,5)+0\cdot s(2,5))\Bigr )
=1,$$ 
where $\chi$ is the odd, primitive and quartic character modulo $5$ given by 
$\chi (1)=1$, 
$\chi (2)=\zeta_4$, 
$\chi (3) =-\zeta_4$ 
and $\chi (4)=-1$ 
($2$ being a generator of the multiplicative group $({\mathbb Z}/5{\mathbb Z})^*$ 
we may chose for $\chi$ the character defined by $\chi(2^k) =\zeta_4^k$).\\
In the same way, for $K$ the imaginary quartic subfield of ${\mathbb Q}(\zeta_{13})$ we have $Q_K=1$, $w_K=2$ 
and since $s(1,13)=11/13$, $s(2,13)=4/13$, $s(3,13)=1/13$, $s(4,13)=-1/13$ and $s(5,13)=0$ we get
$$h_K^-
=\sum_{n=1}^5\Re(\chi (n))s(n,13)
=11/13+1/13+1/13=1
=1,$$ 
where $\chi$ is the odd, primitive and quartic character modulo $13$ such that 
$\chi (1)=1$, 
$\chi (2)=\zeta_4$, 
$\chi (3) =1$ 
and $\chi (4)=-1$ 
($2$ being a generator of the multiplicative group $({\mathbb Z}/13{\mathbb Z})^*$ 
we may chose for $\chi$ the character defined by $\chi(2^k) =\zeta_4^k$).

We came up with [BE] and [GR] while working on \cite{LMQJM} and \cite{LMCJM}, 
where the interested reader will find new results and conjectures on Dedekind sums.

\bibliography{central}

\begin{thebibliography}{??????}
\bibitem[Apo]{Apo} 
M. T. Apostol.
\textit {Modular functions and Dirichlet series in number theory}.
Graduate Texts in Mathematics {\bf 41}. Springer-Verlag, New York, 1976.

\bibitem[BE]{BE}
Bruce C. Berndt and Ronald J. Evans.
\newblock Dedekind sums and class numbers. 
\newblock {\em Monatsh. Math.} {\bf 84} (1977), 265--273. 

\bibitem[GR]{GR}
Larry Joel Goldstein and Michael Razar.
\newblock A generalization of Dirichlet's class number formula. 
\newblock {\em Duke Math. J.} {\bf 43} (1976), 349--358. 

\bibitem[LM21]{LMQJM}
S. Louboutin and M. Munsch.
\newblock Second moment of Dirichlet $L$-functions, 
character sums over subgroups 
and upper bounds on relative class numbers. 
\newblock {\em Quart. J. Math.} {\bf 72} (2021), 1379--1399.

\bibitem[LM2?]{LMCJM}
S. Louboutin and M. Munsch.
\newblock Mean square values of $L$-functions over subgroups for non primitive characters, 
Dedekind sums 
and bounds on relative class numbers.
\newblock {\em Canad. J. Math.}, to appear.

\bibitem[Lou94]{LouCMB36/37}
S. Louboutin.
\newblock Quelques formules exactes pour des moyennes de fonctions $L$ de Dirichlet.
\newblock {\em Canad. Math. Bull.} \textbf{36} (1993), 190--196.
\newblock Addendum.
\newblock {\em Canad. Math. Bull.} \textbf{37} (1994), p. 89.

\bibitem[Lou19]{LouBKMS56}
S. Louboutin.
\newblock On the denominator of Dedekind sums.
\newblock {\em Bull. Korean Math. Soc.} {\bf 56} (2019), 815--827.

\bibitem[RG]{RG}
H. Rademacher and E. Grosswald.
\newblock Dedekind sums.
\newblock {\em The Carus Mathematical Monographs,} {\bf 16}. 
\newblock {\em The Mathematical Association of America, Washington, D.C.,} 1972.

\bibitem[Was]{Was} 
L. C. Washington.
\newblock {\em Introduction to Cyclotomic Fields.} 
\newblock Second Edition. Graduate Texts in Mathematics {\bf 83}. Springer-Verlag, New York, 1997.
\end{thebibliography}

\end{document}